\documentclass[12pt]{amsart}
\usepackage{amssymb,amsmath,epsfig,graphics,latexsym,psfrag}
\usepackage[english]{babel}
\usepackage[all]{xy}
\usepackage{amscd, amssymb, amsmath, amsthm, graphics}
\usepackage{amsmath,amsfonts,amsthm,amssymb}
\usepackage{latexsym,amsmath}
\usepackage{graphicx,psfrag}
\usepackage{mathrsfs}

\newtheorem{theorem}{Theorem}[section]
\newtheorem{proposition}[theorem]{Proposition}
\newtheorem{corollary}[theorem]{Corollary}
\newtheorem{lemma}[theorem]{Lemma}

\theoremstyle{definition}

\theoremstyle{remark}
\newtheorem{remark}[theorem]{Remark}

\newtheorem{question}{Question}
\numberwithin{equation}{section}
\renewcommand{\t}{ \widetilde}
\renewcommand{\hat}{ \widehat}
\renewcommand{\b}{ \partial}
\newcommand{\Z}{\bf Z}
\newcommand{\R}{\bf R}
\newcommand{\N}{\bf N}

\newcommand{\Q}{\bf Q}
\newcommand{\Hi}{\bf H}

\newcommand{\C}{\bf C}

\renewcommand{\S}{\bf S}
\renewcommand{\l}{\langle}
\renewcommand{\r}{\rangle}
\newcommand{\e}{\varepsilon}
\newcommand{\z}[1]{{\Z}/#1{\Z}}
\renewcommand{\o}{\overline}
\newcommand{\abs}[1]{\lvert#1\rvert}

\newcommand{\co}{\colon\thinspace}
\renewcommand{\epsilon}{\varepsilon}
\renewcommand{\c}{\mathcal}

\usepackage{times}

\begin{document}
\sloppy

\title[Graph manifolds have  virtually positive Seifert volume]{Graph manifolds have  virtually positive Seifert volume}

\author{Pierre Derbez}
\address{Centre de Math\'ematiques et d'Informatique,
Technopole de Chateau-Gombert,
39, rue Fr\'ed\'eric Joliot-Curie -
 13453 Marseille Cedex 13}
\curraddr{}
\email{derbez@cmi.univ-mrs.fr}

\author{Shicheng Wang}
\address{Department of Mathematics, Peking University, Beijing, China}
\email{wangsc@math.pku.edu.cn}



\date{\today}

\begin{abstract}
 This paper shows that the Seifert volume of each closed non-trivial graph manifold is virtually positive.
 As a consequence, for each closed orientable prime $3$-manifold  $N$,
 the set of mapping degrees $\c{D}(M,N)$ is finite for any $3$-manifold $M$, unless $N$ is finitely covered by
  either a torus bundle, or a trivial circle bundle, or the $3$-sphere.
\end{abstract}
\maketitle

\vspace{-.5cm} \tableofcontents

\section{Introduction}
Recall that the Thurston's geometrization conjecture, which has been
verified by Perelman, claims that the interior of each component of
the Jaco-Shalen-Johannson decomposition  of a  prime 3-manifold
supports a geometric structure  (i.e. admits a complete locally
homogeneous metric). In other words each JSJ-piece is either
hyperbolic or Seifert.

After a work of Milnor and  Thurston \cite{MT},  Gromov
\cite{G} introduced  the simplicial volume $\|N\|$ of a manifold
$N$, which turns out to be very important. This invariant detects
exactly the hyperbolic part of a
 3-manifold $N$, i.e. it is proportional to the sum of the hyperbolic volume of the
hyperbolic pieces in the JSJ decomposition of $N$, therefore it is
 not affected by the gluing of the JSJ-pieces since it is additive with respect to amenable splittings  (see \cite{G} and also \cite{Th},
\cite{So}). Moreover $||*||$ does respect  the mapping degrees, i.e.
 for any map $f\co M\to N$ then
 $||M||\geq\abs{{\rm deg}(f)}||N||$.

Inspired by \cite{G}, another invariant, the  so-called Seifert
volume $SV(*)$, was introduced by Brooks and Goldman  during their
study of the Godbillon-Vey invariant of foliations \cite{BG1}
\cite{BG2}, which also respects the mapping degrees. It vanishes for closed $3$-manifolds supporting the geometries ${\S}^3, {\S}^2\times{\R}, {\R}^3, {\rm\bf Nil}, {\Hi}^2\times{\R}, {\rm\bf Sol}$  and it is positive
for Seifert manifolds supporting the $\widetilde{\rm SL_2}(\R)$
geometry.

The Seifert volume is rather strange and very little is known about
it. It can be either zero or non-zero for hyperbolic 3-manifolds
\cite{BG2}; and it was not addressed for 3-manifolds with
non-trivial JSJ-decomposition since it was introduced.  It seems
that the Seifert volume reflects the gluing information of the JSJ
pieces.
Call a prime $3$-manifold $N$  a graph manifold, if each JSJ piece
of $N$ is a Seifert manifold. Hence graph manifolds are precisely
those prime 3-manifolds with  zero simplicial volume. Call a graph
manifold  trivial if it is covered by either a torus bundle or a
Seifert manifold. Our main result is

\begin{theorem}\label{main}
For any closed non-trivial graph manifold $N$, there exists a finite
covering $\t{N}$ of $N$ whose Seifert volume $
SV(\t{N})$ is positive.
\end{theorem}

\begin{remark}Theorem \ref{main} is known for graph manifolds which are
homology 3-spheres \cite{Re}.
\end{remark}

Next we will apply this theorem to the study of mapping degrees. Let
$M$ and $N$ be two closed orientable $3$-dimensional manifolds. Let
$\c{D}(M,N)$ be the set of degrees of maps from $M$ to $N$, that is
$$\c{D}(M,N)=\{d\in{\Z}\,| f\co M\to N,\, \, \,{\deg}(f)=d\}.$$

According to \cite{CT}, Gromov think it is a fundamental problem in
topology to determine the set $\c{D}(M,N)$. Indeed the supremum of
absolute values of degrees in $\c{D}(M,N)$ has been addressed by
Milnor and Thurston in 1970's \cite{MT}. A basic property of
$\c{D}(M,N)$ is reflected  in the following
\begin{question}\label{gromov} (see  also \cite[Problem A]{Re} and  \cite[Question 1.3]{W2}):
For which closed orientable $3$-manifolds $N$, is the set
$\c{D}(M,N)$  finite for any given closed orientable $3$-manifold
$M$?
\end{question}

 With the properties of the simplicial volume and
the Seifert volume discussed above, as well as some elementary
constructions, the answer to Question 1 is known for prime
3-manifolds excepted for  non-trivial graph manifolds. It is
expected that $\c{D}(M,N)$ is finite for any $M$ when $N$ is a
non-trivial graph manifold \cite{W2}, \cite{DW}. As an application of Theorem
\ref{main}  we get
\begin{theorem}\label{$D(M,N)$}
For each closed non-trivial graph manifold $N$ the set $\c{D}(M,N)$
of mapping degrees is finite for any $3$-manifold $M$.
\end{theorem}

\begin{remark} Theorem \ref{$D(M,N)$} is known when $M=N$, see \cite{W2},
\cite{D2};  and when $M$ is  graph manifold by using standard form
for maps between graph manifolds, see \cite{DW}, \cite{D1}.
\end{remark}

As a corollary we get a complete answer to Question 1 for prime $3$-manifolds:
\begin{corollary}\label{cor}
For each closed prime $3$-manifold $N$, the set $\c{D}(M,N)$ is finite for
any $M$ if and only if $N$ is not finitely covered by either a torus
bundle, or a trivial circle bundle, or ${\S}^3$.
\end{corollary}

Note that the Seifert volume $SV$ of a closed orientable $3$-manifold
$N$ is defined as the supremum of the volume of representations 
$\phi\co\pi_1N\to{\rm Iso}\t {\rm SL_2}$, where ${\rm Iso}\t {\rm SL_2}$ 
denotes the isometry group of $\t{\rm SL_2}(\R)$. For technical reason, to prove
the non-vanishing property stated in Theorem \ref{main}, we only need
to consider those representations into $\t{{\rm SL}_2}(\R)$, viewed as a
subgroup of ${\rm Iso}\t {\rm SL_2}$. Precisely Theorem \ref{main} is a
consequence of the following

\begin{theorem}\label{connection}
For any closed non-trivial graph manifold $N$, there exists a finite
covering $\t{N}$ of $N$ and a representation $\phi: \pi_1\t{N}\to
\t{{\rm SL}_2}(\R)$  whose volume is nonzero.
\end{theorem}

The paper is organized as follows.

In section \ref{2} we collect the notions and invariants needed in our
proofs: transversely projective foliations, flat connections, volume
of representation, Godbillon-Vey and Chern-Simons invariants. Next
we specialize our discussion to representations  into $\t{{\rm
SL}_2}(\R)$ for which, there are many interpretations of the volume of a
representation, for instance as a Godbillon-Vey invariant of
transversely projective foliations or as a Chern-Simons invariant of
flat connections. These notions have already been introduced and
studied in many important papers, but since they are crucial for our
approach, we recall and organize them for the convenience of the
reader.

In Section \ref{3} we define some technical  objects for graph manifolds:
their coordinates and gluing matrices, canonical framings,  the
absolute Euler number, the separable and  characteristic coverings.
Then we construct some finite coverings of graph manifolds which are
needed in the proof of Theorem \ref{connection}.

Section \ref{4} is devoted to the proof of Theorem \ref{connection} and
the main results. The main difficulty in the proof of Theorem
\ref{connection} is that it seems not easy to know the behavior of the
volumes of representations   with respect to the JSJ-decomposition
of a general graph manifold. Thus our strategy consists to use the
foliation and their Godbillon-Vey invariant to compute the local
volume of a representation (when restricted to the geometric pieces
of the manifold) and next we consider the associated local  flat
connections that one can glue together using gauge transformations
to construct a global representation. One can ensure that the volume
of our representation is nonzero because one can estimate the
Chern-Simons invariant of the associated connection.

For standard terminologies in topology and geometry of 3-manifolds,
see \cite{JS}, \cite{Sc} and \cite{Th}, in foliation see \cite{CC}.
Otherwise in each subsection we list the main references.

From now on, all manifolds are orientable, all coefficients are in
${\R}$, unless otherwise specified. In particular $\t{{\rm SL}_2}(\R)$, ${\rm SL_2}(\R)$
and  ${\rm PSL}_2(\R)$ are denoted by $\t{{\rm SL}_2}$, ${\rm SL}_2$ and ${\rm
PSL}_2$ respectively.

\section{Volumes of representations of  $3$-manifold groups into  $\t{{\rm SL}_2}$}\label{2}

\subsection{Geometry, algebra and topology around Lie group $\widetilde
{\rm {SL}_2}$} The main reference in this subsection is \cite{Sc}\label{2.1}.

 We are going to explain the information
 contained in the following commutative diagram of exact sequences, which will be useful throughout the paper.

\[ \begin{CD}
1 @> >> {\Z} @> >> \widetilde {{\rm SL}_2}@>p >>{\rm PSL}_2@> >> 1\\
@.  @V\cap  VV @V\cap   VV @V \cong  VV\\
1 @> >> {\R} @> >> {\rm Iso}_e \widetilde {{\rm SL}_2} @> \t{p}  >> {\rm PSL}_2@> >>1\\
@.  @A  AA @A  \psi AA @A  \bar\psi AA\\
1 @> >> {\Z} @> >> \pi_1N @> p_*>> \pi_1O_N @> >>1\\
\end{CD} \,\,\,\,\,\,(1) \]

Let ${\rm Iso}_+{\Hi}^2$ be the group of orientation preserving  isometries of
the 2-dimensional hyperbolic space ${\Hi}^2$.  It we identify ${\Hi}^2$ with
the upper-half plane then, we have ${\rm Iso}_+{\Hi}^2={\rm PSL}_2$. This
identification provides the well known action of ${\rm PSL}_2$ on the circle
${\S}^1={\R}\cup \{\infty\}$ (up to conformal transformation on ${\S}^2={\C}\cup\{\infty\}$).

We denote by $p: \t{{\rm SL}_2}\to {\rm PSL}_2$ the universal
 covering, where $p$ is the infinite cyclic covering. Then $\t{{\rm SL}_2}$ is a 3-dimensional 1-connected Lie group,
which provides one of Thurston's eight geometries, which is derived
as below: First ${\rm PSL}_2$, identified with the unit tangent
bundle of ${\Hi}^2$ is automatically geometrized, then $\t{{\rm
SL}_2}$ is endowed with the lifted geometry, which is invariant
under the left multiplication of $\t{{\rm SL}_2}$, which allows to view $\t{{\rm SL}_2}$ as a subgroup of ${\rm Iso}
\widetilde {{\rm SL}_2}$. The circle bundle
structure of the unit tangent bundle of ${\Hi}^2$  provides
 a line bundle structure over ${\Hi}^2$ on $\t{{\rm SL}_2}$, which
 is
often identified with ${\rm Int}{\bf D}^2\times{\R}$

 The action of  ${\rm PSL}_2$
  on the unit circle ${\S}^1$  is lifted to the action of  $\t{{\rm
  SL}_2}$ on ${\R}$. Let $ {\c{R}}_{\alpha}
\in {\rm PSL}_2$ be the rotation of the unit disc of angle
$2\pi\alpha$.  Then $\c{R}_{\alpha}$ admits a lift on the form
$x\mapsto x+2\pi\alpha$ denoted by ${\rm sh}(\alpha)$. The center of
$\t{{\rm SL}_2}$ can be identified with the set $\{{\rm sh}(n),
n\in{\Z}\}$ which is isomorphic to ${\Z}$. The set which consists of
all shift maps $x\mapsto x+\alpha$, when $\alpha$ runs over ${\R}$
will be termed \emph{the core of} $\t{{\rm SL}_2}$ and denoted by
${\rm Core}\left(\t{{\rm SL}_2}\right)$, which corresponds to the
real axis $\{0\}\times{\R}$ in ${\rm Int}{\bf D}^2\times{\R}$.

As we have seen that  $\widetilde {{\rm SL}_2}$ is a subgroup of ${\rm Iso}
\widetilde {{\rm SL}_2}$, the whole ${\rm Iso} \widetilde {{\rm SL}_2}$ is a
4-dimensional Lie group of two components. The identity component ${\rm Iso}_e \widetilde{{\rm SL}_2}$ of ${\rm Iso}
\widetilde {{\rm SL}_2}$ (which preserves the orientation of the lines), is
obtained by extending the action of ${\Z}$, the center of $ \widetilde
{{\rm SL}_2}$, on the lines, to the action of ${\R}$ to the lines. All those
geometric and algebraic informations are surveyed in the first two
lines of table (1).

 Let  $F_{g,n}$ be an
oriented  $n$-punctured surface of genus $g$ with boundary components $s_1,...,s_n$
with $n> 0$. Then $N'=F_{g,n}\times{\S}^1$ is oriented if ${\S}^1$ is
oriented. Let  $h_i$ be the oriented ${\S}^1$ fiber on the torus
$s_i\times h_i$
 (call such pairs $\{(s_i, h_i)\}$ a section-fiber coordinate
system). Let $0\leq l\leq n$. Now attach $l$ solid tori $V_i$ to the
boundary tori of $N'$ such that the meridian of $V_i$ is identified
with the slope $r_i=s_i^{\alpha_i}h_i^{\beta_i}$ where $\alpha_i>0,
(\alpha_i,\beta_i)=1$ for $i=1,...,l$. Denote the resulting manifold
by $\left(
g,n-l;\frac{\beta_1}{\alpha_1},\cdots,\frac{\beta_l}{\alpha_l}\right)$
which has the Seifert fiber structure extended from the circle
bundle structure of $N'$. Each orientable Seifert fibered space with
orientable base $F_{g,n-l}$ and with $l$
exceptional fibers is obtained in such a way.

If $N$ is closed, i.e. if $l=n$, then  define the Euler number of
the Seifert fibration by
$$e(N)=\sum_{i=1}^l \frac{\beta_i}{\alpha_i}\in{\Q}$$
and  the Euler characteristic of the orbifold $O(N)$ by
$$\chi_{O(N)}=2-2g-\sum_{i=1}^l\left(1-\frac 1{\alpha_i}\right)\in{\Q}.$$

A closed orientable 3-manifold $N$ supports the $\widetilde
{{\rm SL}_2}$-geometry, i.e.  there is a discrete faithful representation
$\psi: \pi_1N\to {\rm Iso} \widetilde {{\rm SL}_2}$, if and only if $N$
 is a Seifert manifold with non-zero Euler number $e(N)$ and
negative Euler characteristic $\chi_{O(N)}$. Moreover $\psi$ induces
an injection from the third line to the second line presented in
table (1).

\subsection{Some fact about flat bundles, Chern-Simons invariants}  The main references in this subsection are \cite{Mi},
\cite{DFN}, \cite{KK}, \cite{CC}\label{2.2} and also \cite{BG1}. 
 Let $M$ be a smooth $n$-dimensional manifold and  $G$ be a semi-simple Lie group which acts
 smoothly on some $p$-dimensional manifold $X$.
 Then to each representation $\phi\co\pi_1M\to G$ one can associate a flat fiber  bundle over
 $M$ with fiber $X$ and group $G$ constructed as follows: $\pi_1M$ acts on the universal
 covering $\t{M}$ of $M$ by deck transformations and on $X$ via $\phi$ so it acts diagonally on the product
 $\t{M}\times X$ and we can form the quotient $M\times_{\phi}X=\t{M}\times X/\pi_1M$ which is the flat $X$-bundle over
 $M$ with structure group $G$  corresponding to $\phi$. We define the Euler class of $\phi$, denoted by $\t{e}(\phi)$, to be the Euler
 class $\t{e}$  of  $M\times_{\phi}X\to M$.  This bundle is endowed with a co-dimension
 $p$ foliation $\mathfrak{F}$ transverse to the fiber $X$ whose leaves are the images of
 the leaves $\t{M}\times\{\ast\}$ under the projection. Such a
foliation will be termed \emph{horizontal}. The
 foliation is naturally  endowed with a $(G,X)$-transverse structure, in the sense of \cite[paragraph 1]{BG1}, which means  that the transition
 functions of $\mathfrak{F}$ can be taken to lie in $G$. In other words $G$ is the structure group
 of $\mathfrak{F}$ and $\phi$ is the holonomy of $\mathfrak{F}$.

 When $X=G$ and  acts on
itself by right  multiplication $R_g$ then $M\times_{\phi}G$ is a
principal $G$-bundle over $M$.  Here we obtain a foliation $\mathfrak{F}$ in
$M\times_{\phi}G$ which is invariant under the natural action of $G$
on $M\times_{\phi}G$. By \cite[Lemma 25.1.3]{DFN}, this gives rise to a \emph{flat 
connection} on $M\times_{\phi}G$  defined by an
equation $\Omega=0$ where $\Omega$ is a flat $G$-invariant
differential $1$-form which takes values in the Lie algebra
$\mathfrak{g}$ of $G$.

We are interested in this paper with the trivialized principal $G$-bundle $P=M\times G$ over $M$. We denote by $\c{A}$ the set of  connections on $P$. This corresponds  to $1$-forms $\Omega$ defined on $P$ taking value on $\mathfrak{g}$ such that $R^{\ast}_g\Omega={\rm ad}(g^{-1})\Omega$ for any $g\in G$ and which are the "identity" on the vector tangent to the fiber $G$. The equation $\Omega=0$ gives rise to a horizontal ${\rm codim}(\mathfrak{g})$-plan field $\{H_p\}_{p\in P}$, right invariant under the $G$-action: $H_{p.g}=R_g^{\ast}(H_p)$. We denote by $F^{\Omega}=d\Omega+\Omega\wedge\Omega$ the \emph{curvature of} $\Omega$. It is well known that the equation $F^{\Omega}=0$ is equivalent to the fact that the field $\ker\Omega$ is integrable and we say that the connection is flat. We denote by $\c{A}_0$ the set of flat connections. The set $\c{A}$ is acted by the gauge group $\c{G}$  of all bundle automorphisms covering the identity and this action preserves the subset $\c{A}_0$.    Given an element $\Omega\in\c{A}_0$ one can define a representation $\phi\co\pi_1M\to G$ by "lifting the loops in the leaves". More precisely there is a embedding     
$$\c{A}_0/\c{G}\hookrightarrow{\rm Hom}(\pi_1M\to G)/{\rm Conjugation}$$
that is generally not surjective. Since the bundle $P$ is trivialized one can idenfity  $\c{A}_0$ with the set $\Omega^1(M)\otimes\mathfrak{g}$ of $1$-forms on $M$ taking values in $\mathfrak{g}$  and  we can identify $\c{G}$  with the set ${\rm Maps}(M\to G)$.

When $M$ is an orientable compact $3$-manifold, the Chern-Simons invariant of an element $A\in\Omega^1(M)\otimes\mathfrak{g}$ is defined by the formula
$$\mathfrak{cs}_M(A)=\int_M{\rm tr}\left(dA\wedge A+\frac{2}{3}A\wedge A\wedge A\right)$$

Note that when $G$ is contractible then the correspondance $$\c{A}_0/\c{G}\hookrightarrow{\rm Hom}(\pi_1M\to G)/{\rm Conjugation}$$ is actually a bijection. Indeed given a representation $\phi\co\pi_1M\to G$ then
since $G$ is contractible, the bundle $M\times_{\phi}G$ is trivial and one can choose a bundle isomorphism which is actually a trivialization  $\tau\co M\times_{\phi}G\to P$ which defines a flat connection $A$ on $P$ in the pre-image of $\phi$. Note that this construction does not depend on the choice of $\tau$ because if we choose an other trivialization $\tau'$ defining a connection $B$ on $P$ then $\tau'\tau^{-1}\in\c{G}$ and thus $A$ and $B$ are gauge equivalent.

  On the other hand, when $G$ is contractible and when $M$ is a closed orientable $3$-manifold then   $\mathfrak{cs}_M(A)$ is gauge invariant. Indeed let $A$ and $B$ be two flat connections on $M$ such that $B=gA$ for some $g\in\c{G}$. Since $G$ is contractible, there exists a path $(g_t)_{0\leq t\leq 1}$ of gauge transformations such that $g_0=g$ and $g_1(x)=1$ when $x\in M$ and where $1$ denotes the identity element of $G$. Then the path of flat  connections $A_t=g_t.A$ between $A$ and $B$ defines a connection $\mathbb A$ on $M\times I$   whose curvature is denoted by $F^{\mathbb A}=d{\mathbb A}+{\mathbb A}\wedge{\mathbb A}$. We apply now the arguments in \cite[paragraph 4]{KK}:  since $${\rm tr}\left(F^{\mathbb A}\wedge F^{\mathbb A}\right)=d{\rm tr}\left(d{\mathbb A}\wedge {\mathbb A}+\frac{2}{3}{\mathbb A}\wedge {\mathbb A}\wedge {\mathbb A}\right)$$  then by the Stokes formula
  $$\mathfrak{cs}_M(B)-\mathfrak{cs}_M(A)=\int_{M\times I}{\rm tr}\left(F^{\mathbb A}\wedge F^{\mathbb A}\right)$$
  On the other hand since $F^{A_t}=0$ for any $t$ then $F^{\mathbb A}\wedge dt=0$ and thus $F^{\mathbb A}\wedge F^{\mathbb A}=0$ which shows that 
  $\mathfrak{cs}_M(A)=\mathfrak{cs}_M(B)$.
\subsection{Eisenbud-Hirsch-Neumann's extension of Milnor-Wood's Theorem.}\label{2.3}
 The main references in this subsection are \cite{Mi}, \cite{Wo} and \cite{EHN}.

Recall that ${\rm PSL}_2$ acts on ${\S}^1$. Inspired by Milnor's pioneer work
\cite{Mi}, Wood prove the following theorem in \cite{Wo},

\begin{theorem}\label{Milnor-Wood} Suppose $M$ is an orientable  circle bundle over a closed
surface $F_g$ of genus $g>0$ with  fiber $h$. Then the
following are equivalent:

(1) The bundle is induced by a representation $\phi: \pi_1(F_g)\to
{\rm PSL}_2$;

(2) There is a representation $\t{\phi} : \pi_1(M)\to
\widetilde{{\rm SL}_2}$ such that $\t{ \phi} (h)={\rm sh}(1)$;

(3) There is a $({\rm PSL}_2,{\S}^1)$ horizontal foliation on $M$;

(4) $|e(M)|\le |\chi(F_g)|$.
\end{theorem}

Later Eisenbud-Hirsch-Neumann's made the following extension
\cite{EHN} of Milnor-Wood's Theorem, which will be used in the
paper:

\begin{theorem}\label{Eisenbud-Hirsch-Neumann} Suppose $M$ is a closed orientable Seifert manifold with
a regular fiber $h$ and  base of genus $>0$. Then

(1) (\cite{EHN} Corollary 4.3) There is a $({\rm PSL}_2,{\S}^1)$
horizontal foliation on $M$ if and only if there is a representation
$\t{\phi} : \pi_1(M)\to \widetilde{{\rm SL}_2}$ such that $\t{\phi}
(h)= {\rm sh}(1)$;

(2) (\cite{EHN} Theorem 3.2 and Corollary 4.3) Suppose $M=(g,0; \beta_1/\alpha_1,..., \beta_n/\alpha_n)$, then there is a $({\rm PSL}_2,{\S}^1)$ horizontal
foliation on $M$ if and only if

$$\sum \llcorner{ {\beta_i}/{\alpha_i}}\lrcorner \le -\chi(F_g); \,\,\, \sum \ulcorner{\beta_i}/{\alpha_i}\urcorner \ge
\chi(F_g),$$ where $\llcorner{r }\lrcorner$ and  $\ulcorner
r\urcorner$ for $r\in{\R}$  denote respectively, the greatest
integer $\le r$ and least integer $\ge r$.

\end{theorem}

The proofs of the equivalences of (1), (2) and (3) in Theorem
\ref{Milnor-Wood} and in Theorem \ref{Eisenbud-Hirsch-Neumann} (1)
are conceptual, and the proof the equivalences of (2) and (4) in
Theorem \ref{Milnor-Wood} and in Theorem
\ref{Eisenbud-Hirsch-Neumann} (2) based on an important result:

\begin{proposition} \label{moving distance} (\cite{EHN} Lemma 2.1, Theorem 2.3 and Theorem 4.1)
Let $g_1,...,g_s$ denote elements of $\t{\rm SL_2}$ such that each $g_i$ is conjugated to an element of the form ${\rm sh}(\alpha_i)$, for some $\alpha_i\in \R$. Then the product $g_1...g_s$ can be presented as a product of $g>0$  commutators
$\prod_{i=1}^g [v_i, w_i]$  in $\t{{\rm SL}_2}$ if and only if $$|\alpha_1+...+\alpha_s| < 2g-1$$
\end{proposition}

\subsection{The volume of a representation via transverse
foliations}\label{foliation}

The main references in this subsection are \cite{BG1} and \cite{GV}.

Let $\mathfrak{F}$ be a co-dimensional 1 foliation on a closed
smooth manifold $M$ determined by a 1-form $\omega$. Then by
Frobenius Theorem one have $d\omega= \omega\wedge\delta$ for some
1-form $\delta$. It is observed by Godbillon and Vey \cite{GV}
that the 3-form $\delta \wedge
 d\delta$ is closed and the class $[\delta \wedge
 d\delta]\in H^3(M, \R)$ depends only on the foliation $\mathfrak{F}$,
 which is called the Godbillon-Vey class of $\mathfrak{F}$, denoted by
 $GV(\mathfrak{F})$. It can be verified that

 \begin{lemma}\label{map}(\cite{BG1}) Suppose the map $f: M\to N$ is  transverse to a
 foliation $\mathfrak{F}$ on $N$, then
 $GV(f^*\mathfrak{F})=f^*GV(\mathfrak{F})$, where
 $f^*\mathfrak{F}$ is the induced foliation on $M$. Moreover if $\c{F}$ has a $(G,X)$-transverse
 structure then $f^*\mathfrak{F}$ also has a $(G,X)$-transverse structure.
 \end{lemma}

\begin{lemma}\label{Euler} (Proposition 1 \cite{BG1}) Suppose $\mathfrak{F}$
  is a $({\rm PSL}_2,{\S}^1)$-horizontal foliation
on a flat circle  bundle ${\S}^1\to E \to M$. Then
 $$\int_{\S^1} GV(\mathfrak{F})=4\pi^2\t{e}(E)$$
 where $\int_{\S^1}\co H^3(E)\to H^2(M)$ denotes the integration along the fiber.
 \end{lemma}

Let $M$ be a closed orientable 3-manifold and $\phi\co\pi_1M\to
{\rm PSL}_2$ be a representation with zero Euler class. Since ${\rm PSL}_2$
acts to ${\S}^1$ then one can consider the corresponding flat circle
bundle $M\times_{\phi}{\S}^1$ over $M$ and the associated horizontal
$({\rm PSL}_2,{\S}^1)$-foliation $\mathfrak{F}_{\phi}$. Since the
Euler class of $\phi$ is zero we can choose a section $\delta $ of
$M\times_{\phi}{\S}^1\to M$. Brooks and Goldman (see \cite[Lemma 2]{BG1})
show that $\delta^{\ast}GV(\mathfrak{F}_{\phi})$ only depends on
$\phi$ (and not on a chosen section $\delta$). Then they can define
the volume of a representation  also called the Godbillon Vey
invariant of $\phi$ by setting
$$GV(\phi)=\int_M\delta^{\ast}GV(\mathfrak{F}_{\phi})$$

\begin{remark}
Assume that $M$ is already endowed with a $({\rm PSL}_2,{\S}^1)$-foliation
$\mathfrak{F}$. Then there is a canonical
way to define a flat ${\S}^1$-bundle $E$ over $M$ with structure  group ${\rm PSL}_2$, a
horizontal foliation  $\mathfrak{F}'$ on $E$,  and a section $s:
M\to E$ such that $\mathfrak{F}=s^*\mathfrak{F}'$.
  Then if $\phi$   denotes associated representation   we get (\cite{BG1} Lemma 1)
 $$GV(\phi)=\int_MGV(\mathfrak{F})$$
\end{remark}

\subsection{The volume of a representation via continuous
cohomology}\label{volume} The main references in this subsection are
\cite{BG1}, \cite{BG2} and \cite{V} (for some portions see also
\cite{Re}, \cite{WZ}).

The volume   $GV(\phi)$ can be obtained after integrating a volume form induced
by $\phi$ from the ${\rm Iso}\t{{\rm SL}}_2$-invariant ${\t{\rm
SL}}_2$-volume form. We first state the following Van Est Theorem,
for the definition of continuous cohomology $H_{\text{cont}}^*(G)$
of a Lie group $G$, see \cite{BG1}.

\begin{theorem} \cite{V}
Let $G$ be a connected Lie group, then $H_{\text{cont}}^*(G)$ is
isomorphic to the cohomology of $G$-invariant forms on $G/K$, where
$K$ is a maximal compact subgroup of $G$.
\end{theorem}

Suppose $(G,X)$ is a contractible $3$-dimensional geometry. Then by definition $X=G/G_0$, where $G_0$ is a maximal compact subgroup of $G$. Let
 $\phi:\pi_1(M)\to G$ be a representation, where $M$ is a closed orientable $3$-manifold. Then  there is a natural map
$$\phi^*\co H_{\text{cont}}^*(G)\simeq H^*(\text{$G$-invariant forms on $X$})\to H^*(M)$$ obtained
as follows:    Let $q: \t{
M}\times X\to X$ be the projection, where $\t{
M}$ is the universal cover of $M$. For each $G$-invariant closed
form $\omega$ on $X$, $q^*(\omega)$ is a $\pi_1(M)$-invariant closed
form on $\t{M}\times X$ (for the diagonal action), which induces a form $\omega'$ on
$M\times _\phi X$. Then $s^*(\omega')$ is a closed form on $M$,
where $s: M\to M\times _\phi X$ is a section, and we define
$\phi^*(\omega)=s^*{\omega}'$ (since $X$ is contractible, such a
section exists and  all such sections are homotopic).

Suppose that $(G,X)= ({\rm Iso}\widetilde {{\rm SL}_2}, \widetilde {{\rm SL}_2})$ and  $\omega$ is the unique $G$-invariant
volume 3-form on $X$, then  the volume of the representation $\phi$ is defined by

$${\rm Vol}(\phi)=\int _M \phi^*(\omega)$$

From the definition, the following fact is elementary

\begin{lemma}\label{factorization}
If a  representation $\phi:\pi_1(M)\to G$ factors through
$\pi_1(S)$, then $\phi^* : H^*_{\text{cont}}(G)\to H^*(M)$ factors
through $H^*(S)$.
\end{lemma}

For a given representation $\phi\co\pi_1M \to{\rm PSL}_2$, $\phi$
lifts to $\t\phi\co\pi_1M \to{\t{\rm SL}_2}$ if and only if
$\t{e}(\phi)=0$ in $H^2(M,{\Z})$ ; and
 $\phi$ lifts to $\t\phi\co\pi_1M \to{\rm Iso}\t{{\rm SL}_2}$ if and only
if $\t{e}(\phi)$ is torsion  in $H^2(M,{\Z})$. 
$GV(\phi)$ defined for $\phi\co\pi_1M\to{\rm PSL}_2$
with zero Euler class in \cite{BG1} was extended further by Brooks and
Goldman to those
with torsion Euler class in \cite{BG2}, by proving that
$GV(\mathfrak{F}_{\phi})$ can be considered as a volume form on $M$ as just
discussed.
It is proved that $GV(\phi)$ takes only finitely many valus when
$\phi\co\pi_1M\to {\rm PSL}_2$ runs over the representations with torsion
Euler classes \cite[Theorem 1]{BG2}. Then they define the so-called
Seifert volume of $M$ by setting
$$SV(M)=\max\left\{|GV(\phi)|, Ê \phi\co\pi_1M\to {\rm PSL}_2, \t{e}(\phi)\ {\rm is\ torsion}\right\}$$

We recall the following very
useful fact which has been verified in \cite{BG1} and \cite{BG2}.

\begin{proposition}\label{CV=Vol} Let $M$ be a closed oriented
$3$-manifold, $\phi\co\pi_1M\to {\rm PSL}_2$ be a representation with
zero, resp. torsion, Euler class  and $\t\phi\co\pi_1M \to
 \t {{\rm SL}_2}$, resp., $\t\phi\co\pi_1M \to
 {\rm Iso}\t {{\rm SL}_2}$   is a lift of $\phi$. Then
$$GV(\phi)={\rm Vol}(\t\phi).$$
\end{proposition}

This gives an alternative definition of the Seifert volume by setting 
$$SV(M)=\max\left\{|{\rm Vol}(\phi)|, Ê \phi\co\pi_1M\to {\rm Iso}\t{\rm SL_2} \right\}$$
\begin{theorem}\label{fonct}\cite[Lemma 3]{BG2}
 If $f: M\to N$ is a  map of degree $d \ne 0$ between closed
orientable $3$-manifolds, then $SV(M)\geqslant d SV(N)$.
\end{theorem}

\begin{remark}\label{identification}
For the content below in this paper, we will restrict on
representations into $\t {{\rm SL}_2}$, or equivalently on 
representations into ${{\rm PSL}_2}$ with zero Euler class. To
simplify the notation, when we consider the Godbillon Vey Invariant,  we often do not distinct  a representation into
${\rm PSL}_2$ and its lifts into $\t{{\rm SL}_2}$, or alternatively a
representation into $\t {{\rm SL}_2}$ and its projection into ${{\rm PSL}_2}$.
\end{remark}

\subsection{Relative version of Godbillon-Vey invariants}\label{chern} The main references in this subsection are
\cite{Kh} and \cite{KK}. The following result is proved by Khoi \cite[Section 4]{Kh}.

\begin{proposition}\label{vtk} Let $M$ be a closed oriented $3$-manifold. Let $\phi\co\pi_1M \to{\t{\rm
 SL}_2}$  denote a representation. Then for any flat connection $A\in\Omega^1(M)\otimes\mathfrak{sl}_2$ corresponding to $\phi$  we have
$GV(\phi)=
2\mathfrak{cs}_M(A)$.
\end{proposition}

Note that this identification is also suggested in \cite{Re}.

Let $X$ denote  a connected smooth oriented $3$-manifold with smooth toral
boundary $\b X=T_1\cup...\cup T_r$. Assume each $T_i$ is endowed
with a pair of essential simple closed curves $(\lambda_i,\mu_i)$ which meet
transversally at a single point $\star$. Choose coordinate such that the
exponential map ${\R}^2\to T_i$ sends the horizontal line to
$\lambda_i$ and the vertical line to $\mu_i$. We identify a regular
neighborhood of $T_i$  with the product $T_i\times[0,1]$ such that
$T_i=T_i\times\{1\}$. The $1$-forms $dx$ and $dy$ on ${\R}^2$
descend to $1$-forms on $T_i$ that we  denote still by $dx$ and
$dy$. This allows to define an oriented basis of $1$-forms near
$T_i$ by $\{dx,dy,dr\}$ where $r\in[0,1]$. The purpose of this
section is to state the construction of the \emph{relative Godbillon-Vey invariant}
of an  \emph{elliptic representation} $\rho$ of $\pi_1X$ into
$\t{{\rm SL}_2}$. This
construction naturally extends the construction of the closed case.

Pick some base point $x_0$ in ${\rm int}(X)$ and consider a representation
$\rho\co\pi_1(X,x_0)\to\t{{\rm SL}_2}$ and suppose that for each
$i=1,...,r$, there exists $-\infty<\alpha_i,\beta_i<+\infty$ such
that $\rho(\lambda_{i,x_0})$ is conjugated to ${\rm sh}(\alpha_i)$ and
$\rho(\mu_{i,x_0})$ is conjugated to ${\rm sh}(\beta_i)$, where $\lambda_{i,x_0}$ and
$\mu_{i,x_0}$ are the elements of $\pi_1(X,x_0)$ corresponding to
$\lambda_i$ and $\mu_i$ respectively (more precisely  take a path from $x_0$ to $\lambda_i\cap\mu_i=\{\star\}$ to view these elements in $\pi_1(X,x_0)$). We say $\rho$ is an
\emph{elliptic representation}.

Note that in the proof of Lemma \ref{vtk}, Khoi uses the  standard identification of ${\rm SL}_2({\R})$ with the
group ${\rm SU}(1,1)$ (by conjugation in ${\rm SL}_2({\C})$) whose
Lie Algebra is given by
$$\left\{\begin{pmatrix}
i\alpha&\beta \\
\o{\beta}&-i\alpha
\end{pmatrix}, \alpha\in{\R},\beta\in{\C}\right\}.$$
In the sequel we will keep the same presentation for convenience.  Denote by $A_{\rho}$ a connection on the trivialized bundle $X\times\t{\rm SL_2}$ corresponding to $\rho$. Since $\rho$ is elliptic this implies, using the computation of \cite[paragraph 3 (1)]{Kh} and standard gauge theory that for each $i=1,...,r$ there exists a gauge transformation $g_i\co T_i\times[0,1]\to\t{\rm SL_2}$ such that 
$g_i.A_{\rho}|T_i\times[0,1]$ is equal to $$\begin{pmatrix}
i\alpha_i&0 \\
0&-i\alpha_i
\end{pmatrix}\otimes dx+\begin{pmatrix}
i\beta_i&0 \\
0&-i\beta_i
\end{pmatrix}\otimes dy+0\otimes dr$$
Since $\t{\rm SL_2}$ is contractible then one can extend the map $g_1\coprod ...\coprod g_r\co\b X\times[0,1]\to\t{\rm SL_2}$ to $X$. 
Thus using  \cite[Theorem 4.2]{Kh} we get:

\begin{lemma}\label{nf}
Let $\rho$ be an elliptic representation of $\pi_1X$   and let $A_{\rho}\in\Omega^1(X)\otimes\mathfrak{sl}_2$ be a flat connection corresponding to $\rho$. Then

(i) there exists a gauge transformation $g\co X\to\t{\rm SL_2}$ such that
$$J^{\ast}_igA_{\rho}= \begin{pmatrix}
i\alpha_i&0 \\
0&-i\alpha_i
\end{pmatrix}\otimes dx+\begin{pmatrix}
i\beta_i&0 \\
0&-i\beta_i
\end{pmatrix}\otimes dy+0\otimes dr$$ where  $J_i\co T_i\times[0,1]\to X$ denotes the inclusion. In this case we say that the connection is in \emph{elliptic normal  form}.

(ii) If $A$ and $B$ denote two connections corresponding to $\rho$
in elliptic normal forms then
$\mathfrak{cs}_X(A)=\mathfrak{cs}_X(B)$.
\end{lemma}
As suggested in  \cite{Kh}, we define   the Godbillon-Vey invariant
of an elliptic representation of $\pi_1X$ into $\t{{\rm SL}_2}$,
denoted by $GV(\rho)$, by setting $GV(\rho)=2\mathfrak{cs}_X(A)$
where $A$ is any flat connection in elliptic normal form 
representing $\rho$. This definition makes sense by Lemma \ref{nf}.
The relation between Godbillon-Vey invariant and Chern-Simons
invariant is very convenient for our purpose since $\mathfrak{cs}$
is (by definition) additive with respect to a given decomposition of
a manifold.

\section{Coordinated graph manifolds and their finite coverings}\label{3}

\subsection{Coordinated graph manifolds}\label{3.1} A reference for this subsection can be found in \cite{DW}.

 Let $N$ be an oriented  non-trivial graph manifold. Denote by
$\c{T}_N$ the family of JSJ tori of $N$, by $N^{\ast}$ the set
$N\setminus \c{T}_N=\{\Sigma_1,...,\Sigma_n\}$ of the JSJ pieces of
$N$, by $\tau\co\b N^{\ast}\to\b N^{\ast}$ the associated sewing
involution defined in \cite{JS}.

Since Theorem \ref{connection}  is stated "up to finite coverings"
there are no loss of generality assuming that $N$ satisfies the
following properties (*) (for details see \cite{DW}):

(i) each Seifert piece $\Sigma_i$ of $N$ is an orientable  ${\S}^1$-bundle over
an orientable surface $F_i$ of genus at least $2$,

(ii) each characteristic torus of $N$ is shared by two distinct
Seifert pieces of $N$.

A {\it dual graph}  of $N$, denoted by $\Gamma_N$, is given as
follows: each vertex represents a JSJ piece of $N$; each edge
represents a JSJ torus of $N$; an edge $e$ connects two vertices
$v_1$ and $v_2$  if and only if the corresponding
JSJ torus is shared by the corresponding JSJ pieces.

Call a dual graph $\Gamma_N$  {\it directed} if each edge of
$\Gamma_N$ is directed, in other words, is endowed with an arrow.
Once $\Gamma_N$ is directed, the sewing involution $\tau$ becomes a
well defined map, still denoted by $\tau\co\b N^{\ast}\to\b
N^{\ast}$.

Note that each $\Sigma_i$ admits a unique Seifert fibration, up to
isotopy. Denote by $h_i$ the homotopy class of the regular fiber of
$\Sigma_i$. Since $N$ is non-trivial then $\b\Sigma_i\not=\emptyset$ and thus  there exists a cross section $s_i\co F_i\to\Sigma_i$. Call
$\Sigma_i$  {\it coordinated}, if

(1) such a section $s_i\co F_i\to\Sigma_i$ is chosen,

(2) both $F_i$ and $h_i$ are oriented so that their product
orientation is matched by the orientation of $\Sigma_i$ induced by
that of $N$.

Once $\Sigma_i$ is coordinated, then  $\partial F_i$ and the  fiber $h_i$
gives a basis of $H_1(T;{\Z})$ for each component $T$ of $\partial
\Sigma_i$. We also say that $\Sigma_i$ is endowed with a $(s,h)$-basis.

Call $N$ is {\it coordinated}, if each component $\Sigma_i$ of $N^*$
is coordinated and $\Gamma_N$ is directed.

Once $N$ is coordinated, then each torus $T$ in $\c{T}_N$ is
associated with a unique $2\times 2$-matrix $A_T$ provided by the
gluing map $\tau| \co T_-(s_-, h_-)\to T_+(s_+,h_+)$: where $T_-,
T_+$ are two torus components in $\partial N^*$ provided by $T$,
with basis $(s_-, h_-)$ and $(s_+,h_+)$ respectively, and
$$\tau (s_-, h_-)= (s_+,h_+)A_T.$$

Call $\{A_T, T\in \c{T}\}$ the {\it gluing matrices}.

Let $\Sigma$ denote an orientable Seifert manifold with regular
fiber $h$. A \emph{framing} $\alpha$ of $\Sigma$ is to assign a
simple closed essential curve not homotopic to the regular fiber of
$\Sigma$, for each component $T$ of $\b\Sigma$. Denote by
$\Sigma(\alpha)$ the closed Seifert 3-manifold obtained from
$\Sigma$ after Dehn fillings along the framing  $\alpha$.   Each Seifert piece $\Sigma_i$ of $N^{\ast}$ is
endowed with a \emph{canonical framing} $\alpha_i$ given by
the regular fiber of the Seifert pieces of $N^{\ast}$ adjacent to
$\Sigma_i$. Denote by $\hat{\Sigma}_i$ the space
$\Sigma_i(\alpha_{i})$. By minimality of JSJ decomposition,
$\hat{\Sigma}_i$ admits a unique Seifert fibration extending that of
$\Sigma_i$. One can define the {\it absolute Euler number} of $N$ by
setting
$$|e|(N)=\sum_{i=1}^n|e(\hat{\Sigma}_i)|$$

   Let
$\c{T}$ be a union of tori and let $m$ be a positive integer. Call a
covering $p: \t{\c{T}}\to \c{T}$  {\it $m$-characteristic} if for
each component $T$ of $\c{T}$ and for each component $\t{T}$ of
$\t{\c{T}}$ over $T$,  the restriction $p| : \t{T}\to T$  is the
covering map associated to the characteristic subgroup of index $m
\times m$ in $\pi_1T$. Call a finite covering $\t{N}\to N$ of a
graph manifold $N$ $m$-characteristic if its restriction to
$\c{T}_{\t{N}}\to\c{T}_N$ is {\it $m$-characteristic}.

Next we define the \emph{separable coverings}. Let $\Sigma$ be a
component of $N^{\ast}$ with a given section $s\co F\to\Sigma$, then
$\Sigma\simeq F\times{\S}^1 $.  Let $p\co\t{\Sigma}\to\Sigma$ denote a
finite covering. Since $p$ is a fiber preserving map one can define
the \emph{vertical degree} of $p$ as the integer $d_v$ such that
$p_{\ast}(\t{h})=h^{d_v}$, where $h$ and $\t{h}$ denote the homotopy
class of the regular fiber in $\Sigma$ and $\t{\Sigma}$, and the
\emph{horizontal degree} $d_h$ as the degree of the induced covering
$\bar p: \t{F}\to F$, where $\t F$ denotes the base of the bundle
$\t{\Sigma}$. We have ${\rm deg}(p)=d_v\times d_h$.

On the other hand, $p$ induces a finite covering
$p|_b\co\t{\c{F}}\to F$, with $\t{\c{F}}$ connected, corresponding
to the subgroup $s^{-1}_{\ast}p_{\ast}(\pi_1\t{\Sigma})$.  If ${\rm
deg}(p|_b)= d_h$, then we say that the covering $p$ is
\emph{separable}.

From now on we assume the graph manifold $N$ is coordinated. Let
$p\co\t{N}\to N$ be a finite covering of graph manifolds. Then
obviously $\Gamma_{\t{N}}$ can be directly in a unique way such that
the induced map $p_\#: \Gamma_{\t{N}}\to \Gamma_N$ preserves the
directions of the edges. Below we also assume that $\Gamma_{\t{N}}$
is directed in such a way.

Let $p\co\t{N}\to N$ be a finite covering of graph manifolds. Call
$p$  {\it separable} if the restriction $p| \co\t{\Sigma}\to\Sigma$
on connected Seifert pieces is separable for all possible
$\t{\Sigma}$ and $\Sigma$. Call a coordinate on $\t{N}$  a {\it
lift} of the coordinate of $N$, if for each possible covering $p|
\co\t{\Sigma}\to\Sigma$ on connected Seifert pieces, the
$(s,h)$-basis of $\t{\Sigma}$ is obtained as the pre-image of the $(s,h)$-basis of $\Sigma$.

The following three lemmas are proved in \cite{DW}:

\begin{lemma}\label{separable}
Let  $p\co\t{\Sigma}\to\Sigma$  be a finite covering between
oriented Seifert manifolds.

(i) If $p$ has fiber degree one, then $p$ is a separable covering.

(ii)  If $\Sigma=F\times{\S}^1$ and
 $p$ is a regular covering corresponding to an epimorphism
$\phi\co\pi_1\Sigma=\pi_1F\times{\Z}\to G=G_1\times G_2$ satisfying
$\phi(\pi_1F\times\{1\})=G_1$ and $\phi(\{1\}\times{\Z})=G_2$ then
$p$ is separable.
\end{lemma}

\begin{lemma}\label{stable} Le $N$ be a graph manifold satisfying property (*).

(i) If $p\co\t{N}\to N$ is a separable finite covering then the
coordinate of $N$ can be lifted on $\t{N}$.

(ii) Moreover, if the covering $p$ is characteristic, then for each
component $T$ of $\c{T}_N$ and for each component $\t{T}$ over $T$
we have $A_T=A_{\t{T}}$, where the coordinate of $\t{N}$ is
lifted from $N$.
\end{lemma}

\begin{lemma}\label{unif}
Let $N$ be a closed graph manifold satisfying property (*). For any finite covering $\t{N}\to N$ then $|e|(N)=0$ iff $|e|(\t{N})=0$.  If
$|e|(N)=0$ then $N$ admits a finite covering $\t{N}$ which can be
coordinated such that each gluing matrix is in the form $\pm J$ where 
$$J=\begin{pmatrix}
0&1 \\
1&0
\end{pmatrix}$$
\end{lemma}

\subsection{Constructions of finite coverings}\label{3.2}
Let $N$ be a graph manifold with property (*).

\begin{lemma}\label{genus} Let $S$ be a Seifert piece of $N$ endowed with a framing $\mu_S$.
 For any integer $n\geq 1$ there exists a  separable and characteristic finite covering $p_n\co N_n\to N$ satisfying the following properties :

(i) for each Seifert piece $S'$ (resp. each canonical torus $T'$) of
$N$, the space $p_n^{-1}(S')$, resp. $p_n^{-1}(T')$) is connected,

(ii)
$e\left(p_n^{-1}(S)(p_n^{-1}(\mu_S)\right)=e({S}(\mu_S))$  when $n$ runs over ${\N}$,

(iii) $\left|\chi\left(F_n\right)\right|>n$, where $F_n$ denotes the
base of the Seifert bundle $p_n^{-1}(S)$.
\end{lemma}
\begin{proof}
 After fixing a section of the
bundle $S$ we may identify $S$ with a direct product
$F\times{\S}^1$, where $F$ is an orientable surface with
 genus $\geq 2$. Then $\pi_1F$ admits a presentation
$$\left\l a_1,b_1,...,a_g,b_g,d_1,...,d_p | [a_1,b_1]....[a_g,b_g]d_1...d_p=1\right\r$$ if $F$ has $p\geq 1$ boundary components. Choose a prime number $q>p$ and consider the homomorphism
$$\e_1\co\pi_1F\to\z{q}$$ defined by $\e_1(a_i)=\e_1(b_i)=\o{0}$ for any $i=1,...,g$ and
$\e_1(d_1)=...=\e_1(d_{p-1})=\o{1}$ and the obvious homomorphism
$$\e_2\co\pi_1{\S}^1\to\z{q}$$ Thus define $$\varphi=\e_1\oplus\e_2\co\pi_1S\to\z{q}\oplus\z{q}$$
It follows from our construction that $\e_1(d_p)=-\o{p-1}$.  Since
$q$ and $p-1$ are coprime numbers then $\e_1(d_i)$ is of order $q$
for any $i=1,...,p$. Then the  covering $p_q\co S_q\to S$
corresponding to $\varphi$  is $q$-caracteristic and
separable by Lemma \ref{separable}. Then it follows directly that
for each component $T$ of $\b S$ then $p_q^{-1}(T)$ is connected. We perform this construction for each Seifert piece  of $N$.

Since the  covering $p_q\co S_q\to S$ is $q$-characteristic
and separable for each Seifert piece $S$ of $N$, we can define a
$q^2$-fold covering still denoted by $p_q\co N_q\to N$ by gluing the
coverings $S_q$ together using a lifting of the sewing involution of
$N$ by Lemma \ref{stable} (i). Note that $p_q$ satisfies point (i)
of Lemma \ref{genus} and it is separable and characteristic by our
construction. Moreover  the Euler number
$e({S_q(p_q^{-1}(\mu_S)}))=e({S}(\mu_S))$ does not depend on $q$ by
Lemma \ref{stable} (ii). This proves Lemma \ref{genus} (ii). The
genus $g_q$ of the base of the Seifert bundle $S_q$ is given by the
Riemann-Hurwitz formula
$$2(g_q-g)= (2g+p-2)(q-1)$$
where $g$ denotes  the genus of the base of the Seifert bundle $S$.
Since $g\geq 2$, $p\geq 1$ and $q$ is an arbitrarily large prime
number, this completes the proof of Lemma \ref{genus}.

\end{proof}
\begin{lemma}\label{genusbis}
Let $\Sigma$ be a Seifert piece of $N$ and denote by $S_1,...,S_p$
the adjacent Seifert pieces of $N$. Then for any integer $n$ there
exists a finite separable $1$-characteristic covering $p_n\co N_n\to N$
inducing the trivial covering over $\Sigma\cup\c{T}_N$ and such that
for each $i$, $p_n^{-1}(S_i)$ is a connected circle bundle over a
surface of genus $>n$.
\end{lemma}
\begin{proof}  After fixing a section of the bundle $S_i$ we may identify $S_i$
with a direct product $F_i\times{\S}^1$, where $F_i$ is a surface of genus $g\geq 2$.  Choose a prime number $q$ and consider the homomorphism
$$\e_1\co\pi_1F_i\to\z{q}$$ defined by $\e_1(a_i)=\e_1(b_j)=\o{0}$ for any $i=2,...,g$ and $j=1,...,g$,
$\e_1(d)=\o{0}$ where $d$ runs over the components of $\b F_i\cap S_i$ and $\e_1(a_1)=\o{1}$. Thus
define $$\varphi=\e_1\times 0\co\pi_1S_i\to\z{q}$$ Since the induced
covering is trivial on the boundary (i.e. $1$-characteristic) then one can glue together one component of each
$S_i$ and $q$ copies of each Seifert piece disctinct from
$S_1,...,S_p$ to define a covering $p_q\co N_q\to N$ satisfying the
first statement of the lemma, using Lemma \ref{separable}. Fix a
Seifert piece $S_i$. The genus $g_q$ of the base of $p_q^{-1}(S_i)$
is given by the Riemann Hurwitz formula
$$g_q-g\geq (g-1)(q-1)$$
 Since $g\geq 2$ and since $q$ is an arbitrarily large prime
number, this completes the proof of the lemma.
\end{proof}

\section{Proof of the main results}\label{4}

With the relations among volume of representations, Godbillon-Vey
invariant of foliations and  Chern-Simons invariants of flat
connections given in Subsections \ref{foliation}, \ref{volume} and
\ref{chern}, (also recall Remark \ref{identification}), to prove
Theorem \ref{connection}, we only need to construct a $\t
{{\rm SL}_2}$-flat connection on certain finite covering of $N$ with
non-vanishing Chern-Simons  invariant.

\subsection{Construction of local flat connection with non-zero Chern-Simons invariants}\label{4.1}
Let $N$ be a graph manifold satisfying property (*) of paragraph 3.1. Assume $N$ is
endowed with a family of curves $(\lambda,\mu)$ such that for each
component $T$ of $\c{T}_N$ the set
$T\cap\{(\lambda,\mu)\}=\{\lambda_T,\mu_T\}$ is a pair of essential simple
closed curves meeting transversally at one point. If $\t{N}\to N$ is
a characteristic separable covering then if $(\t{\lambda},\t{\mu})$ is the pre-image of  $(\lambda,\mu)$ in $\t{N}$ then  for each
component $\t{T}$ of $\c{T}_{\t{N}}$ the set
$\{\t{\lambda}_{\t{T}},\t{\mu}_{\t{T}}\}$ is also a pair of simple
closed curves meeting transversally at one point. In the same way, suppose we are given  a family of curves in $N$ denoted by $c(p,q)$ such that for each component $T$ of $\c{T}_N$ the set $T\cap c(p,q)$ is the slope $p_T\lambda_T-q_T\mu_T$, where $p_T,q_T$ are coprime integers. Then the curves $c(p,q)$ have as a pre-image  a family of curves $\t{c}(p,q)$ in $\t{N}$  such that for each
component $\t{T}$ of $\c{T}_{\t{N}}$ each component of 
$\t{c}(p,q)\cap\t{T}$ is the curve $p_T\t{\lambda}_{\t{T}}-q_T{\t{\mu}}_{\t{T}}$.  When $N$ is
endowed with such a system of curves we denote it by
$(N,\lambda,\mu)$ and we assume each component $T$ of $\c{T}_N$ is
parametrized so that the $1$-form $dx$, resp. $dy$, corresponds to
$\lambda$, resp. $\mu$ and we adopt the same notation for any finite
covering $(\t{N},\t{\lambda},\t{\mu})$ of $(N,\lambda,\mu)$ (this
makes sense since $dx$ and $dy$ comes from the universal covering
${\R}^2$ of $T$). The main purpose of this section is to prove Theorem
\ref{connection}. First of all we will need the following

\begin{proposition}\label{elliptic}
Let $(N,\lambda,\mu)$ be a graph manifold satisfying property (*) endowed with a family of curves $c(p,q)$.
Let $\Sigma$ be a Seifert piece of $N$ such that $c(p,q)\cap\b\Sigma$  define a framing of $\Sigma$. 
Then there exists a separable characteristic finite covering $\t{N}$
of $N$ such that each component $\t{\Sigma}$ over $\Sigma$ admits  a
flat connection $A_{\t{\Sigma}}$ in elliptic normal form such that
$$J^{\ast}A_{\t{\Sigma}}=\begin{pmatrix}
{i}{\alpha_T}&0 \\
0&-{i}{\alpha_T}
\end{pmatrix}\otimes dx+\begin{pmatrix}
{i}{\beta_T}&0 \\
0&-{i}{\beta_T}
\end{pmatrix}\otimes dy$$
 for each component $\t{T}$ over $T$, where $\alpha_T,\beta_T\in{\R}$ satisfy $\alpha_Tp_T=\beta_Tq_T$,
 $J\co\t{T}\times[0,1]\to\t{\Sigma}$ denotes the inclusion and moreover
$$\mathfrak{cs}_{\t{\Sigma}}\left(A_{\t{\Sigma}}\right)=2\pi^2e(\t{\Sigma}(\t{c}(p,q)))$$
where $\t{\Sigma}(\t{c}(p,q))$ denotes the Seifert manifold obtained from $\t{\Sigma}$ after Dehn filling along the slopes $\t{c}(p,q)\cap\b\t{\Sigma}$.\end{proposition}

\begin{proof}
 Note that $\Sigma$ can be
identified with $F\times{\S}^1$ for some compact surface $F$. Applying Lemma
\ref{genus} for $n$ sufficiently large, there exists a finite
covering $\t{N}$ of $N$ such that the component
$\t{\Sigma}$ over $\Sigma$
satisfies the condition
$$\left|\chi\left(\t{F}\right)\right|\gg\left|e\left(\t{\Sigma}(\t{c}(p,q))\right)\right|=\left|e\left({\Sigma}(c(p,q))\right)\right|$$
 where $\t{F}$ denotes the underlying surface of $\t{\Sigma}(\t{c}(p,q))$. 

Let $x_0$ denote a base point in the interior of $\t{\Sigma}$. By
Theorem \ref{Eisenbud-Hirsch-Neumann} (2), $\t{\Sigma}(\t{c}(p,q))$
admits a $({\rm PSL}_2,{\S}^1)$-horizontal foliation $\mathfrak{F}$ and equivalently by Theorem
\ref{Eisenbud-Hirsch-Neumann} (1)  there exists a representation
$$\hat \Phi\co\pi_1(\t{\Sigma}(\t{c}(p,q)),x_0)\to\t{{\rm SL}_2}$$  such that $\hat{\Phi}(h_{x_0})={\rm sh}(1)$,
where $h_{x_0}$ denotes the fiber of the bundle $\t{\Sigma}$, $\hat{\Phi}(p_T\t{\lambda}_{\t{T},x_0}-q_T{\t{\mu}}_{\t{T},x_0})={\rm sh}(0)$ and $\hat{\Phi}(\t{\lambda}_{\t{T},x_0})$ is conjugated to ${\rm sh}(\alpha_T)$, for some $\alpha_T\in\R$ and $\hat{\Phi}(\t{\mu}_{\t{T},x_0})$ is conjugated to ${\rm sh}(\beta_T)$, for some $\beta_T\in\R$. Then we have necessarily $\alpha_Tp_T=\beta_Tq_T$.

Write
$\t{\Sigma}(\t{c}(p,q))$ as $\t{\Sigma}\cup\c{V}$ where $\c{V}$ denotes
the disjoint union of the solid tori corresponding to the Dehn
fillings along the slopes $\t{c}(p,q)\cap\b\t{\Sigma}$. Denote by $\hat{A}$ a flat connection over $\t{\Sigma}(\t{c}(p,q))$ 
corresponding to $\hat{\Phi}$. Then by the same arguments as in paragraph 2.6 there exists a gauge transformation $g\co\b\c{V}\times[0,1]\to\t{\rm SL}_2$, where $\b\c{V}\times[0,1]$ is identified with a regular neighborhood of $\b\t{\Sigma}=\b\c{V}$ in $\t{\Sigma}(\t{c}(p,q))$ with $\b\c{V}\times\{1/2\}=\b\t{\Sigma}$ such that, with respect to
the basis $(\t{\lambda},\t{\mu})$ we get, for each component $\t{T}$ of
$\b\t{\Sigma}$:
 $$g\hat{A}|\t{T}\times[0,1]=\begin{pmatrix}
{i}{\alpha_T}&0 \\
0&-{i}{\alpha_T}
\end{pmatrix}\otimes dx+\begin{pmatrix}
{i}{\beta_T}&0 \\
0&-{i}{\beta_T}
\end{pmatrix}\otimes dy$$   
 Then since $\t{\rm SL}_2$ is contractible one can extend $g$ to a gauge transformation $\hat{g}\co\t{\Sigma}(\t{c}(p,q))\to\t{\rm SL}_2$. Denote by
$A_{\t{\Sigma}}\cup A_0$ the decomposition of $\hat{g}.\hat{A}$ corresponding to the decomposition of $\t{\Sigma}(\t{c}(p,q))$ into the union $\t{\Sigma}\cup\c{V}$  where $A_{\t{\Sigma}}$ is the restriction of $\hat{g}.\hat{A}$  over $\t{\Sigma}$  and $A_0$ is  the restriction of $\hat{g}.\hat{A}$  over $\c{V}$. Note that if follows from our construction that both  $A_{\t{\Sigma}}$ and  $A_0$ are in elliptic normal form.  Then, taking care  that $\c{V}$ is glued to $\t{\Sigma}$ by an orientation reversing homeomorphism, i.e:
 $\t{\Sigma}(\t{c}(p,q))=\t{\Sigma}\cup-\c{V}$,  we have the
equalities:
$$GV(\hat{\Phi})=2\mathfrak{cs}_{\t{\Sigma}(\t{c}(p,q))}(\hat{A})=2\mathfrak{cs}_{\t{\Sigma}}(A_{\t{\Sigma}})-
2\mathfrak{cs}_{\c{V}}(A_0)$$ and by \cite[Proposition 6.1]{Kh},
since $A_0$ is in normal form by our construction, $\mathfrak{cs}_{\c{V}}(A_0)=0$.

It remains to compute $GV(\hat{\Phi})$. Since $\t{\Sigma}(\t{c}(p,q))$
admits a $\t{{\rm SL}_2}$-geometry with  base a $2$-orbifold of
genus at least $2$ then there exists a  fiber degree one finite
covering $p\co\hat{\Sigma}(\t{c}(p,q))\to\t{\Sigma}(\t{c}(p,q))$ such that
$\hat{\Sigma}(\t{c}(p,q))$ is a  circle bundle over a smooth surface
$\hat{F}$. Notice  that since $p$ is a fiber degree one covering,
the bundle  $\hat{\Sigma}(\t{c}(p,q))$ is endowed with a representation
$$\hat{\Psi}=\hat{\Phi}|\pi_1\left(\hat{\Sigma}(\t{c}(p,q))\right)$$ such that $\hat{\Psi}(\hat{h})={\rm sh}(1)$ where
$\hat{h}$ denotes the fiber of $\hat{\Sigma}(\t{c}(p,q))$. Then the
pull back $\hat{\mathfrak{F}}=p^{\ast}(\mathfrak{F})$ is a
horizontal $({\rm PSL}_2,{\S}^1)$ foliation on $\hat{\Sigma}(\t{c}(p,q))$
corresponding to $\hat{\Psi}$. By Lemma \ref{map}, we have
$$GV(\hat{\Psi})=\int_{\hat{\Sigma}(\t{c}(p,q))}G-V(\hat{\mathfrak{F}})={\rm deg}
(p)\int_{\t{\Sigma}(\t{c}(p,q))}G-V({\mathfrak{F}})={\rm
deg}(p)GV(\hat{\Phi})$$ Since $\hat{\mathfrak{F}}$ is a horizontal
$({\rm PSL}_2,{\S}^1)$ foliation on $\hat{\Sigma}(\t{c}(p,q))$,  a circle bundle over
a surface $\hat F$, then  using Lemma \ref{Euler} we get
$$\int_{\hat{\Sigma}(\t{c}(p,q))}G-V(\hat{\mathfrak{F}})=\int_{\hat{F}} \int_{{\S}^1}
G-V(\hat{\mathfrak{F}})=\int_{\hat{F}}4\pi^2\t
e\left(\hat{\Sigma}(\t{c}(p,q))\right)
=4\pi^2e\left(\hat{\Sigma}(\t{c}(p,q))\right)$$ Then by the naturality
property of (rationnal) Euler number of Seifert manifolds  we get,
dividing by the degree of the covering
$$GV(\hat{\Phi})=\int_{\t{\Sigma}(\t{c}(p,q))}G-V({\mathfrak{F}})=4\pi^2e\left(\t{\Sigma}(\t{c}(p,q))\right)$$
Hence
$$\mathfrak{cs}_{\t{\Sigma}}\left(A_{\t{\Sigma}}\right)=\frac 12 GV(\hat{\Phi})=2\pi^2e(\t{\Sigma}(\t{c}(p,q)))$$

 This completes
the proof of the proposition.
\end{proof}

\subsection{Extending the connection from local to global}\label{4.2}
\begin{proof}[Proof of Theorem \ref{connection}]
We may assume that $N$ satisfies property (*). The proof of the
proposition depends on the absolute Euler number of the manifold.

\emph{First Case: $|e|(N)\not=0$}. Thus it follows from the
definition of the absolute Euler number that there exists at least
one Seifert piece, say $\Sigma$, of $N$ such that
$e(\hat{\Sigma})\not=0$, where $\hat{\Sigma}$ denotes the closed
Seifert bundle obtained from $\Sigma$ after Dehn filling along the
fibers of the adjacent Seifert pieces of $N$. Passing to a finite
covering, we may apply Proposition \ref{elliptic} to $\Sigma$ whose each
boundary component $T$ is endowed with the couple of curves $(\lambda,\mu)$
given by the $(s,h)$-basis of the adjacent Seifert pieces and with the curves $c(p,q)=\mu$ so that $p_T=0$ and $q_T=1$ for any $T$ in $\b\Sigma$. Fix a
base point $x_0$ in the interior of $\Sigma$. This gives an elliptic
normal form flat connection $A$ on $\Sigma$ such that
$$J^{\ast}A=\begin{pmatrix}
{i}{\alpha_T}&0 \\
0&-{i}{\alpha_T}
\end{pmatrix}\otimes dx$$
where $J\co
T\times[0,1]\to\Sigma$ is the natural inclusion. 
 \begin{center}
\psfrag{sigma}{$\Sigma$} \psfrag{x}{$x_0$} \psfrag{y}{$y_0$}
\psfrag{t}{$t$} \psfrag{d1}{$d_1$}
\psfrag{d2}{$d_2$} \psfrag{h1}{$h_1$}
\psfrag{h2}{$h_2$} \psfrag{S}{$S$} \psfrag{S'}{$S'$} \psfrag{F}{$F$} \psfrag{sm}{${\rm sewing\ map}$}
\includegraphics{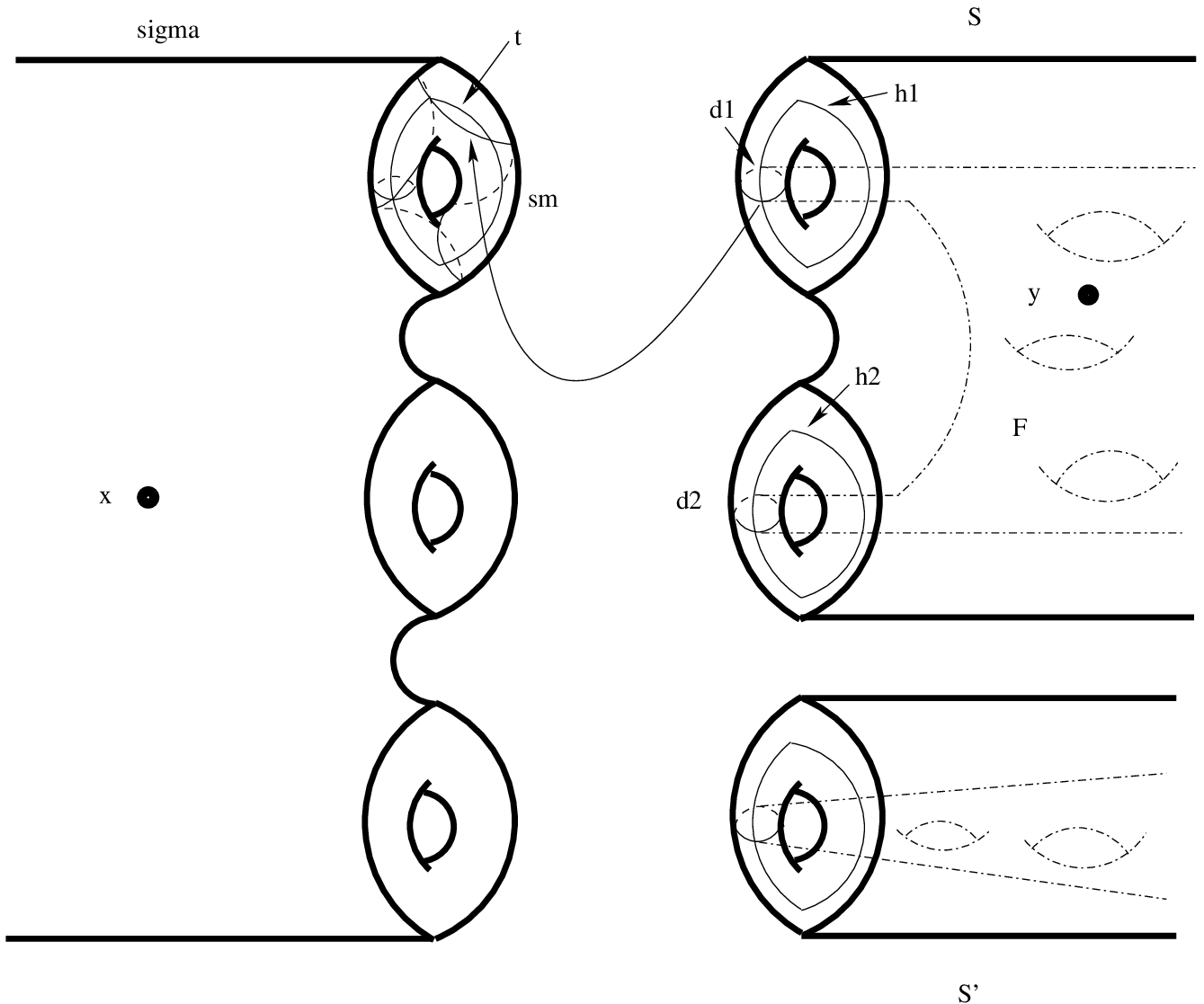}
\centerline{Figure 1}
\end{center}

Denote by $S$ a Seifert bundle of $N$  adjacent to $\Sigma$ with $\b
S\cap\b\Sigma=T_1\cup...\cup T_{r}$ and fix a base point $y_0$ in the
interior of $S$. By property (*), $S$ can be identified with a
product $F\times{\S}^1$. We precise the notation of the above
paragraph : each component $T_i$ viewed as a component of $\b S$ is
endowed with its natural pair of curves $(d_i,h_i)=(\lambda,\mu)\cap T_i$ where
$d_i=T_i\cap\b F$ and $h_i$ is the fiber of $S$ represented in
$T_i$.

 Using Lemma
\ref{genusbis} we may assume that the genus $g(F)$ of each $F$ is
"sufficiently" large. Indeed the coverings of Lemma \ref{genusbis}
are trivial on each components over $\Sigma$ so that one can keep the same flat
connection $A$ and moreover the $(s,h)$ basis can be lifted in such
a way that the gluing matrices do not change (see Figure 2).

\begin{center}
\psfrag{sigma}{$\Sigma$} \psfrag{covering}{$3-{\rm fold\ covering\ of\ } \Sigma\cup S$}  \psfrag{S}{$S$} 
\includegraphics{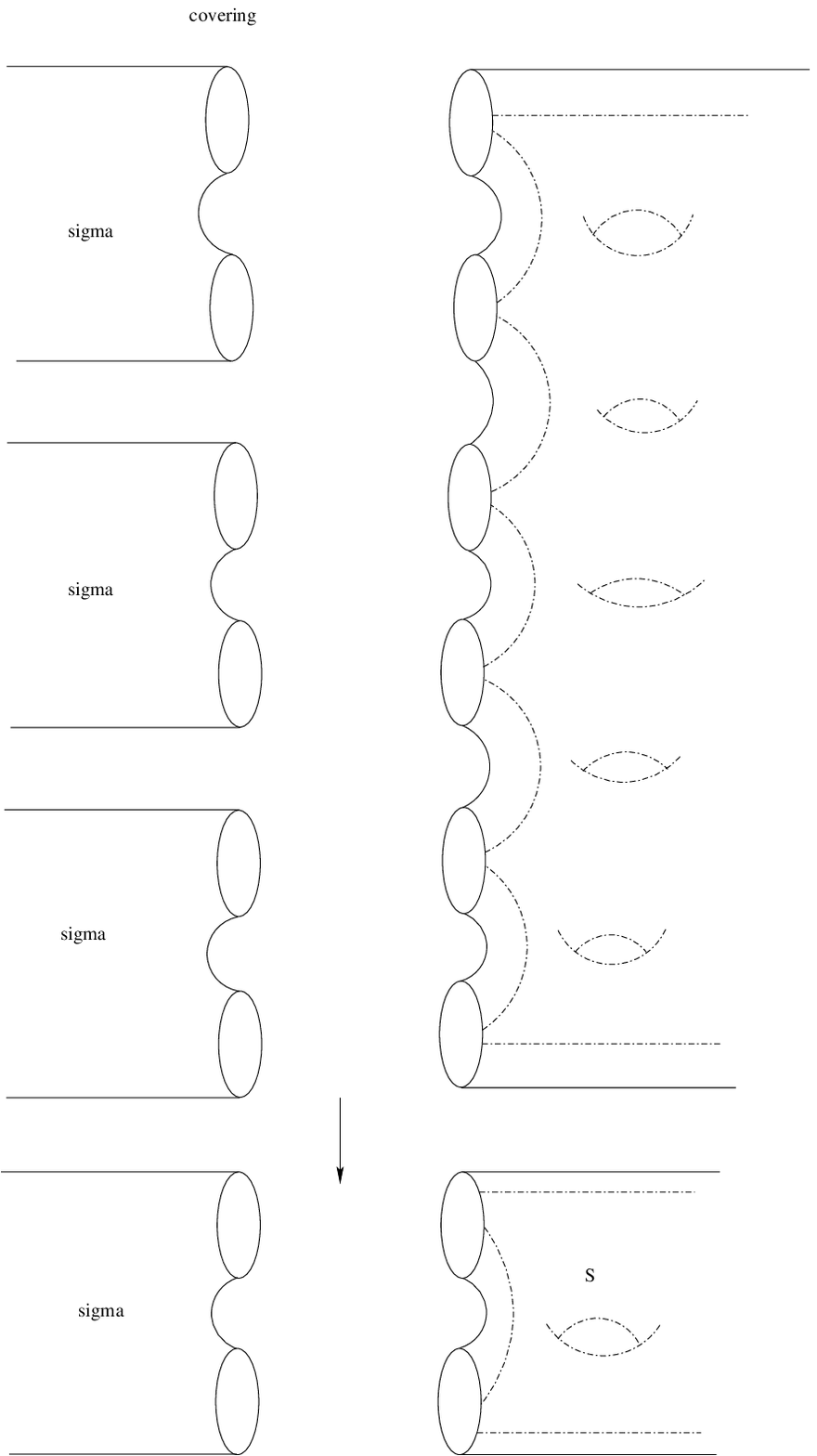}
\centerline{Figure 2}
\end{center}

  Now
applying Proposition \ref{moving distance} we know that
${\rm sh}(\alpha_{T_1}+....+\alpha_{T_r})$ may be written as a product of $g(F)$
commutators $[v_1,w_1]...[v_{g(F)},w_{g(F)}]$. Consider the quotient
of the group $\pi_1(S,y_0)$ given by
$$G=\left\l a_1,b_1,...,a_{g(F)},b_{g(F)},d_{1,y_0},...,d_{r,y_0}|[a_1,b_1]...[a_{g(F)},b_{g(F)}]=d_{1,y_0}...d_{r,y_0}\right\r\times\l h_{y_0}\r$$
obtained by killing the normal subgroup generated by the components
of $\b F$ distinct from $F\cap(T_1\cup...\cup T_r)$. Then the
homomorphism $\Phi_S\co G\to\t{{\rm SL}_2}$ given by
$\Phi_S(h_{y_0})={\rm sh}(0)$, $\Phi_S(d_{i,y_0})={\rm sh}(\alpha_{T_i})$,
$\Phi_S(a_i)=v_i$ and $\Phi_S(b_i)=w_i$ for $i=1,...,g(F)$ is well
defined. Denote by $B$ the flat connection in elliptic normal form
given by Lemma \ref{nf} defined over $S$ corresponding to the
representation $\pi_1S\to G\to \t{{\rm SL}_2}$. More precisely
$J^{\ast}B=0$ when $T$ is distinct from $T_1\cup ....\cup T_r$ and
$$J^{\ast}B=\begin{pmatrix}
{i}{\alpha_T}&0 \\
0&-{i}{\alpha_T}
\end{pmatrix}\otimes dx$$ when $T$ is a component of
$T_1\cup ....\cup T_r$ and $J\co T\times[0,1]\to S$ denotes the
inclusion. Thus we may define a flat connection over $\Sigma\cup S$
denoted by $A\cup B$. Notice that  by a similar
argument as in the proof of Proposition \ref{elliptic}, we have
$$\mathfrak{cs}_S(B)=\mathfrak{cs}_{S(h)}(\hat{B})=\frac
12 GV(\hat\Phi_S),$$ where $\hat{B}$ is obvious extension of $B$ to
$S(h)$ and $\hat\Phi_S: \pi_1(S(h))\to \t{{\rm SL}_2}$ is induced
from $\Phi_S$ since $\Phi_S (h_{y_0})={\rm sh}(0)$. Then  it is easy to check
that we have the factorization
$$\hat\Phi_S: \pi_1(S(h)) \to \pi_1(F)\to \t{{\rm SL}_2}$$
Then by Lemma \ref{factorization}, we have the factorization
$$\hat\Phi_S^*
: H^3_{\text{cont}}(\t{{\rm SL}_2})\to H^3(F)\to H^3(S(h)).$$ Since
$H^3(F)$ vanishes as $F$ is a 2-complex, we must have
$GV(\hat\Phi_S)=0$, see paragraph \ref{volume},
 thus $\mathfrak{cs}_S(B)=\mathfrak{cs}_{S(h)}(\hat{B})=0$.

 We perform this
construction for each Seifert piece of $N$ adjacent to $\Sigma$ and
next we extend the connection by the trivial connection on
$$N\setminus \left(\Sigma\bigcup_{S\in N^{\ast}, \b
S\cap\b\Sigma\not=\emptyset}S\right)$$ This makes sense since the
restrictions of the  representations we constructed on each torus
boundary component $T$ of $S$, where $S\cap \Sigma\ne \emptyset$ and
$T \cap\Sigma= \emptyset$, is trivial. This define a flat connection
on $N$ whose Chern-Simons invariant is equal to
$\mathfrak{cs}_{{\Sigma}}(A)$ which is nonzero by Proposition
\ref{elliptic}. By Lemma \ref{vtk} this completes the proof of the
proposition when $|e|(N)\not=0$.

\vskip 0.5 true cm

\emph{Second Case: $|e|(N)=0$}. Then by Lemma \ref{unif} we may
assume, possibly passing to a finite covering, that the gluing
matrices of $N$ are of the form $\pm\begin{pmatrix}
0&1 \\
1&0
\end{pmatrix}$ with respect to the usual $(s,h)$-basis. Choose two adjacent Seifert pieces in $N$
that we denote by $\Sigma_1,\Sigma_2$ and denote by
$\Sigma=\Sigma_1\cup\Sigma_2$ the connected graph submanifold of $N$
whose Seifert pieces are $\Sigma_1$ and $\Sigma_2$ and by
$\hat{\Sigma}$ the graph manifold obtained from $\Sigma$ after Dehn
fillings along the fiber of the Seifert pieces of $N$ adjacent to
$\Sigma$. 

\begin{center}
\psfrag{n-}{$N\setminus\Sigma$} \psfrag{s1}{$\Sigma_1$}  \psfrag{s2}{$\Sigma_2$} \psfrag{s}{${\rm section\ } s$}
\psfrag{f}{${\rm fiber\ } h$} \psfrag{d1}{$d_1$}
\psfrag{d2}{$d_2$} \psfrag{h1}{$h_1$}
\psfrag{h2}{$h_2$} \psfrag{S}{$S$} \psfrag{F}{$F$} \psfrag{sm}{${\rm sewing\ map}$}
\includegraphics{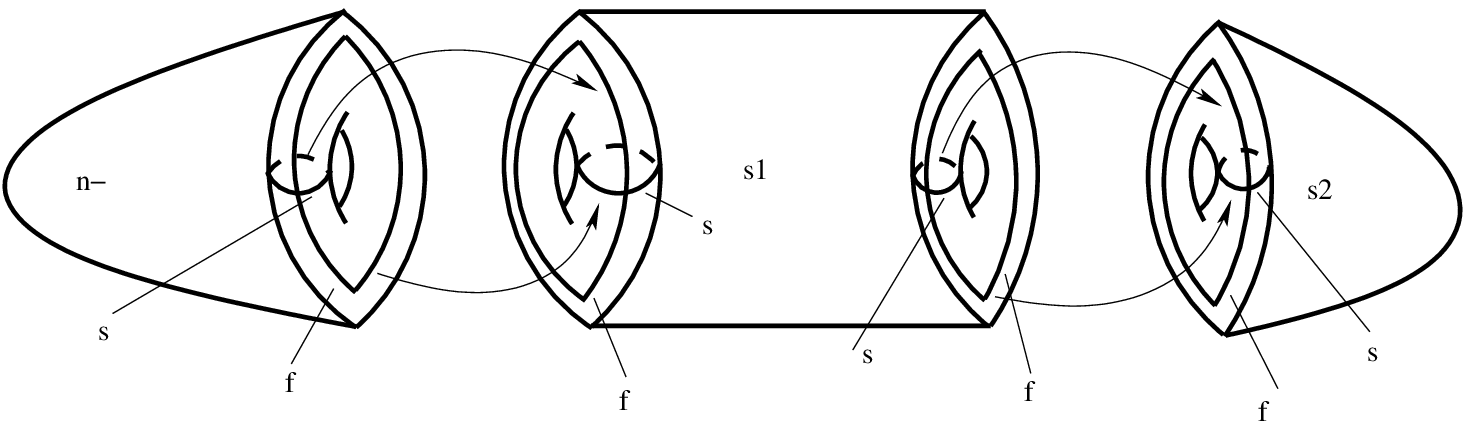}
\centerline{Figure 3}
\end{center}

Denote by $\hat{\Sigma}^{\ast}_1, \hat{\Sigma}^{\ast}_2$ the components of $\hat{\Sigma}^{\ast}=\hat{\Sigma}\setminus\c{T}_{\hat{\Sigma}}$.
For each component $T$ of $\b\hat{\Sigma}^{\ast}_i$, we set $(\lambda,\mu)=(s,h)$ and $c(p,q)=\lambda-\mu$ so that $p_T=q_T=1$ and for  each component $T$  of  $\b\Sigma_i\cap\b{\Sigma}$ then $(\lambda,\mu)$ is given by the $(s,h)$-basis of the  corresponding adjacent
 Seifert pieces of $N\setminus\Sigma$ (see Figure 3) and $c(p,q)=\mu$ so that $p_T=0$ and $q_T=1$. In particular if $\b\hat{\Sigma}^{\ast}_1=\b\hat{\Sigma}^{\ast}_2$ has $r$ components then $$e(\hat{\Sigma}^{\ast}_1(c(p,q)))=-e(\hat{\Sigma}^{\ast}_2(c(p,q)))=\pm r\ \ \ \ \ \ \ (I)$$
The latter equality follows from the special form of gluing
matrices, taking care the orientations.

 Then using Proposition \ref{elliptic} we may assume that for each $i=1,2$ there exists a  flat connection $A_{i}$
  over $\Sigma_i$ in normal elliptic form such that
$$J^{\ast}A_{i}=\begin{pmatrix}
{i}{\alpha_T}&0 \\
0&-{i}{\alpha_T}
\end{pmatrix}\otimes dx+\begin{pmatrix}
{i}{\alpha_T}&0 \\
0&-{i}{\alpha_T}
\end{pmatrix}\otimes dy \ \ \ \ \ \ (**)$$
for each component $T$ of $\b\hat{\Sigma}^{\ast}_i$ where $\alpha_T\in{\R}$,
 $J\co{T}\times[0,1]\to{\Sigma}_i$ denotes the inclusion and 
 $$J^{\ast}A_{i}=\begin{pmatrix}
{i}{\alpha_{T}}&0 \\
0&-{i}{\alpha_{T}}
\end{pmatrix}\otimes dx\ \ \ \ \ \ (***)$$
 for each component  $T$ of $\b\Sigma_i\cap\b{\Sigma}$, where $\alpha_{T}\in{\R}$,
 $J\co{T}\times[0,1]\to{\Sigma}_i$ denotes the inclusion.

  Moreover
$\mathfrak{cs}_{{\Sigma}_i}\left(A_{i}\right)=\pm 2\pi^2 r$. On the
other hand equation (**) allows to define a flat connection
$A=A_1\cup A_2$ over $\Sigma$ such that (by $(I)$)
$$\mathfrak{cs}_{\Sigma}(A)=\mathfrak{cs}_{\Sigma_1}(A_1)-\mathfrak{cs}_{\Sigma_2}(A_2)=\pm 4\pi^2r\not=0$$
As in the \emph{First Case}, equation (***) allows to define a flat
connection $B$ over $N\setminus\Sigma$ in normal elliptic form whose
Chern-Simons invariant vanishes and which can be glued with $A$ on a
regular neighborhood of $\b\Sigma$ defining a flat connection $A\cup
B$ over $N$ such that $\mathfrak{cs}_N(A\cup
B)=\mathfrak{cs}_{\Sigma}(A)\not=0$. This completes the proof
Theorem \ref{connection}.
\end{proof}
\begin{proof}[Proof  of Theorem \ref{main}]
The proof of Theorem \ref{main} follows directly from Theorem
\ref{connection} by the definition of the Seifert volume.
\end{proof}

\begin{proof}[Proof of Theorem \ref{$D(M,N)$}] Let $N$ be a
closed non-trivial graph manifold. By Theorem \ref{main} there is a
finite covering $\t N$ of $N$ with positive Seifert volume $SV(\t
N)$. By Theorem \ref{fonct}, $\c{D}(M,\t N)$ is finite for any
closed $3$-manifold $M$. Then $\c{D}(M,N)$ is finite for any closed
$3$-manifold $M$ by the following

\begin{lemma}\label{covering}\cite[Lemma 3.1]{DW}
Let $p\co N'\to N$ be a finite covering of a closed oriented
$3$-manifold $N$. If $\c{D}(P,N')$ is finite for any closed
$3$-manifold $P$,  then $\c{D}(M,N)$ is finite for any closed
$3$-manifold $M$.
\end{lemma}
\end{proof}
\begin{proof}[Proof of Corollary \ref{cor}]
If $N$ has a positive simplicial volume, the set $\c{D}(M,N)$ is finite for any $M$ since $\|M\|\geq |{\rm deg}(f)|\|N\|$. If $\|N\|=0$ then by Theorem \ref{$D(M,N)$} it remains to consider the case where $N$ is geometric. If $N$ admits the geometry $\t{{\rm SL}_2}$ then $SV(N)>0$ and the set $\c{D}(M,N)$ is finite for any $M$ since $SV(M)\geq |{\rm deg}(f)|SV(N)$. If $N$ admits the geometry, ${\S}^3$ (i.e. $N$ is fnitely covered by ${\S}^3$),
${\S}^2\times{\R}$, ${\Hi}^2\times{\R}$ (i.e. $N$ is finitely covered by a trivial circle bundle)  ${\rm Nil}$, ${\R}^3$,  or
${\rm Sol}$ (i.e. $N$ is finitely covered by a torus bundle) then $\c{D}(N,N)$ is infinite by \cite{W2}. This completes the proof of the corollary.
\end{proof}

\end{document}